\documentclass[12pt]{article}%
\usepackage{amssymb}
\usepackage{amsfonts}
\usepackage{amsmath}
\usepackage{graphicx}%
\providecommand{\U}[1]{\protect\rule{.1in}{.1in}}
\newtheorem{theorem}{Theorem}[section]

\newtheorem{conjecture}[theorem]{Conjecture}
\newtheorem{corollary}[theorem]{Corollary}

\newtheorem{Fact}[theorem]{Fact}
\newtheorem{lemma}[theorem]{Lemma}

\newtheorem{problem}[theorem]{Problem}
\newtheorem{question}[theorem]{Question}
\newtheorem{proposition}[theorem]{Proposition}

\providecommand{\boksie}{\ensuremath{\mathbin{\raisebox{0.3mm}{$\scriptstyle\square$}}}}
\begin{document}

\title{\textbf{Domination, Eternal Domination, and }\\\textbf{Clique Covering}}
\author{William F. Klostermeyer\\University of North Florida\\Jacksonville, FL 32224-2669\\{\small wkloster@unf.edu}
\and C. M. Mynhardt\thanks{Supported by the Natural Sciences and Engineering
Research Council of Canada.}\\Department of Mathematics and Statistics\\University of Victoria, P.O. Box 1700 STN CSC\\Victoria, BC, \textsc{Canada} V8W 2Y2\\{\small kieka@uvic.ca}}
\date{}
\maketitle

\begin{abstract}
Eternal and $\mathrm{m}$-eternal domination are concerned with using mobile
guards to protect a graph against infinite sequences of attacks at vertices.
Eternal domination allows one guard to move per attack, whereas more than one
guard may move per attack in the $\mathrm{m}$-eternal domination model.
Inequality chains consisting of the domination, eternal domination,
$\mathrm{m}$-eternal domination, independence, and clique covering numbers of
graph are explored in this paper.

Among other results, we characterize bipartite and triangle-free graphs with
domination and eternal domination numbers equal to two, trees with equal
$\mathrm{m}$-eternal domination and clique covering numbers, and two classes
of graphs with equal domination, eternal domination and clique covering numbers.

\end{abstract}

\noindent\textbf{Keywords:\hspace{0.1in}}dominating set, eternal dominating
set, independent set, clique cover

\noindent\textbf{AMS Subject Classification Number 2000:\hspace{0.1in}}05C69

\section{Introduction}

\label{Sec_Models}A \emph{dominating set }of a finite, undirected graph
$G=(V,E)$ is a set $D\subseteq V$ such that each vertex in $V-D$ is adjacent
to at least one vertex in $D$. The minimum cardinality amongst all dominating
sets of $G$ is the \emph{domination number}, $\gamma(G)$. By imposing
conditions on the subgraph $G[D]$ of $G$ induced by $D$, one can obtain
several varieties of dominating sets and their associated parameters. For
example, if $G[D]$ is connected, then $D$ is a \emph{connected dominating set
}and the corresponding parameter is the \emph{connected} \emph{domination
number} $\gamma_{c}(G)$.

Domination theory can be considered the precursor to the study of graph
protection: one may view a dominating set as an immobile set of guards
protecting a graph. A thorough survey of domination theory can be found in
\cite{HHS}. In this paper, we consider two forms of dynamic domination which
aim to protect a graph against an infinite sequence of attacks occurring at
the vertices of the graph.

Let $\{D_{i}\}$, $D_{i}\subseteq V$, $i\geq1$, be a collection of sets of
vertices of the same cardinality, with one guard located on each vertex of
$D_{i}$. The two problems considered in this paper can each be modeled as a
two-player game between a \emph{defender} and an \emph{attacker}: the defender
chooses $D_{1}$ as well as each $D_{i}$, $i>1$, while the attacker chooses the
infinite sequence of vertices corresponding to the locations of the attacks
$r_{1},r_{2},\ldots$. Players alternate turns, with the defender first
choosing the initial location of guards. The attacker goes next and chooses a
vertex to attack. Each attack is dealt with by the defender by choosing the
next $D_{i}$ subject to some constraints that depend on the particular game
(see below). The defender wins the game if they can successfully defend any
sequence of attacks, subject to the constraints of the game described below;
the attacker wins otherwise.

We say that a vertex is \emph{protected} if there is a guard on the vertex or
on an adjacent vertex. A vertex $v$ is \emph{occupied} if there is a guard on
$v$, otherwise $v$ is \emph{unoccupied}. An attack at an unoccupied vertex $x$
is \emph{defended} if a guard moves to the attacked vertex. If the guard moves
to $x$ from $v$, we also say $v$ \emph{defends} $x$.

For the \textbf{eternal domination problem}, each $D_{i}$, $i\geq1$, is
required to be a dominating set, $r_{i}\in V$ (assume without loss of
generality $r_{i}\notin D_{i}$), and $D_{i+1}$ is obtained from $D_{i}$ by
moving one guard to $r_{i}$ from an adjacent vertex $v\in D_{i}$. If the
defender can win the game with the sets $\{D_{i}\}$, then each $D_{i}$ is an
\emph{eternal dominating set (EDS)}. The size of a smallest EDS of $G$ is the
\emph{eternal domination number} $\gamma^{\infty}(G)$. This problem was first
studied by Burger et al.~in \cite{BCG2} and will sometimes be referred to as
the \emph{one-guard moves} model. It has been subsequently studied in
\cite{ABB,GK, KM} and other papers.

For the \textbf{m-eternal dominating set problem}, each $D_{i}$, $i\geq1$, is
required to be a dominating set, $r_{i}\in V$ (assume without loss of
generality $r_{i}\notin D_{i}$), and $D_{i+1}$ is obtained from $D_{i}$ by
moving guards to neighboring vertices. That is, each guard in $D_{i}$ may move
to an adjacent vertex, as long as one guard moves to $r_{i}$. Thus it is
required that $r_{i}\in D_{i+1}$. The size of a smallest $\mathrm{m}%
$\emph{-eternal dominating set (}$\mathrm{m}$\emph{-EDS)} (defined similarly
to an EDS) of $G$ is the $\mathrm{m}$\emph{-eternal domination number}
$\gamma_{\mathrm{m}}^{\infty}(G)$. This \textquotedblleft multiple guards
move\textquotedblright\ version of the problem was introduced by Goddard,
Hedetniemi and Hedetniemi \cite{GHH}. We refer to this as the
\textquotedblleft all-guards move" model of eternal domination. This problem
has been subsequently studied in \cite{GK2, KM2} and other papers.

It is clear from the definitions that $\gamma^{\infty}(G)\geq\gamma
_{\mathrm{m}}^{\infty}(G)\geq\gamma(G)$ for all graphs $G$. A survey on
several variations of eternal dominating sets, including the two just defined,
can be found in \cite{KMSurvey}. Our focus in this paper is comparing these
graph protection parameters to other parameters which will be defined and
reviewed in the next section. We pay special attention to the study of graph
classes that satisfy equality in bounds on $\gamma^{\infty}$ and
$\gamma_{\mathrm{m}}^{\infty}$. After providing definitions, background and
known results in Section \ref{SecBack}, we consider $\mathrm{m}$%
\emph{-}eternal domination in graphs with $\alpha=3$ in Section
\ref{Sec_Small_alpha} as initiation of the study of graphs $G$ for which
$\gamma_{\mathrm{m}}^{\infty}(G)=\alpha(G)$. In Section \ref{Sec2} we
characterize bipartite graphs with $\gamma=\gamma^{\infty}$, and bipartite and
triangle-free graphs with $\gamma=\gamma_{\mathrm{m}}^{\infty}=2$. As the main
result of this paper, trees with equal $\mathrm{m}$-eternal domination and
clique covering numbers are characterized in Section \ref{SecTrees}, and in
Section \ref{SecTheta} we consider the problem of whether $\gamma
(G)=\gamma^{\infty}(G)$ implies that $\gamma(G)=\theta(G)$. We end with a
number of open problems and questions in Section \ref{Sec_Open}.

\section{Definitions and Background}

\label{SecBack}The \emph{open} and \emph{closed neighborhoods} of $X\subseteq
V$ are $N(X)=\{v\in V:v$ is adjacent to a vertex in $X\}$ and $N[X]=N(X)\cup
X$, respectively, and $N(\{v\})$ and $N[\{v\}]$ are abbreviated, as usual, to
$N(v)$ and $N[v]$. The set $\overline{N[v]}$ is the set of all vertices not
dominated by $v$. For any $v\in X$, the \emph{private neighborhood}
$\operatorname{pn}(v,X)$\emph{ of }$v$\emph{ with respect to }$X$ is the set
of all vertices in $N[v]$ that are not contained in the closed neighborhood of
any other vertex in $X$, i.e., $\operatorname{pn}(v,X)=N[v]-N[X-\{v\}]$. The
elements of $\operatorname{pn}(v,X)$ are the \emph{private neighbors of }$v$
\emph{relative to }$X$. The \emph{external private neighborhood},
$\operatorname{epn}(v,X)$, is defined similarly, except that $N(v)$ replaces
$N[v]$ in the definition.

In a tree $T$, a \emph{leaf} is a degree one vertex, a \emph{stem} is a vertex
adjacent to a leaf, and a branch vertex is a vertex of degree at least three.
For any $v\in V(T)$, a $v$-\emph{endpath} is a path from $v$ to a leaf, all of
whose internal vertices have degree two in $T$. An \emph{end-branch-vertex} is
a branch vertex $v$ such that exactly one edge incident with $v$ does not lie
on a $v$-endpath. Every tree with at least two branch vertices has at least
two end-branch vertices. A (non-trivial) \emph{star} is a tree $K_{1,r}$,
$r\geq1$.

We denote the minimum and maximum degree of a graph $G$ by $\delta(G)$ and
$\Delta(G)$ respectively, and its independence number by $\alpha(G)$. The
\emph{clique covering number }$\theta(G)$ is the minimum number $k$ of sets in
a partition $V=V_{1}\cup\cdots\cup V_{k}$ of $V$ such that each $G[V_{i}]$ is
complete. Hence $\theta(G)$ equals the chromatic number $\chi(\overline{G})$
of the complement $\overline{G}$ of $G$. Since $\chi(G)=\omega(G)$ (the size
of a maximum clique) if $G$ is perfect, and $G$ is perfect if and only if
$\overline{G}$ is perfect, $\alpha(G)=\theta(G)$ for all perfect graphs.

As first observed by Burger et al.~\cite{BCG2}, $\gamma^{\infty}$ lies between
the independence and clique covering numbers, giving the inequality chain below.

\begin{Fact}
\label{FactED_Bound}For any graph $G$, $\gamma(G) \leq\alpha(G)\leq
\gamma^{\infty}(G)\leq\theta(G)$.
\end{Fact}

Since $\alpha(G)=\theta(G)$ for perfect graphs, the rightmost two bounds in
Fact \ref{FactED_Bound} are tight for perfect graphs. A topic that has
received much attention is finding classes of non-perfect graphs that satisfy
equality in one or more of the bounds in Fact \ref{FactED_Bound}. A number of
graphs classes have been shown to satisfy $\gamma^{\infty}(G)=\theta(G)$, such
as circular-arc graphs \cite{regan} and series-parallel graphs \cite{ABB}. It
is, as of yet, not known whether $\gamma^{\infty}(G)=\theta(G)$ for all planar
graphs $G$.


The following upper bound is due to Klostermeyer and MacGillivray \cite{KM};
Goldwasser and Klostermeyer \cite{GK} show that the bound is sharp.

\begin{theorem}
\label{ThmKM}\emph{\cite{KM}\hspace{0.1in}}For any graph $G$,%
\[
\gamma^{\infty}(G)\leq\binom{\alpha(G)+1}{2}.
\]

\end{theorem}

Goddard et al.~\cite{GHH} determine $\gamma_{\mathrm{m}}^{\infty}(G)$ exactly
for complete graphs, paths, cycles, and complete bipartite graphs. Further,
they show that $\gamma_{\mathrm{m}}^{\infty}(G)=\gamma(G)$ for all Cayley
graphs $G$ obtainable from abelian groups. Their assertion that this equality
holds for all Cayley graphs is shown to be false in \cite{BdSL}.

The inherent symmetry of Cayley graphs provides a sort of foothold for
$\mathrm{m}$-eternal domination; an open problem is to determine other classes
of graphs where $\gamma_{\mathrm{m}}^{\infty}(G)=\gamma(G)$. Goddard et
al.~also prove the following fundamental bound.

\begin{theorem}
\label{GHH1} \emph{\cite{GHH}} For all graphs $G$, $\gamma(G)\leq
\gamma_{\mathrm{m}}^{\infty}(G)\leq\alpha(G)$.
\end{theorem}

In order to get a better upper bound on $\gamma_{\mathrm{m}}^{\infty}$,
Goddard et al.~define a \emph{neo-colonization} to be a partition
$\mathcal{P}=\{V_{1},V_{2},\ldots,V_{t}\}$ of $G$ such that each $V_{i}$
induces a connected graph \cite{GHH}. A part $V_{i}$ is assigned \emph{weight}
$w(V_{i})=1$ if $V_{i}$ induces a clique, and $w(V_{i})=1+\gamma_{c}%
(G[V_{i}])$ otherwise, where $\gamma_{c}(G[V_{i}])$ is the connected
domination number of the subgraph induced by $V_{i}$. The \emph{weight}
$w(\mathcal{P})$ of a neo-colonization $\mathcal{P}$ is the sum of the weights
of its parts. Define $\theta_{c}(G)$ to be the minimum weight of any
neo-colonization of $G$. Goddard et al.~\cite{GHH} prove that $\gamma
_{\mathrm{m}}^{\infty}(G)\leq\theta_{c}(G)\leq\gamma_{c}(G)+1$. In general,
however, $\alpha(G)$ and $\theta_{c}(G)$ are not comparable: consider
$\theta_{c}(K_{1,5})<\alpha(K_{1,5})$, $\theta_{c}(K_{n})=\alpha(K_{n})$, and
$\theta_{c}(C_{5})=3>\alpha(C_{5})=2$. On the other hand, $\theta_{c}%
(G)\leq\alpha(G)$ for all perfect graphs $G$ because $\theta_{c}(G)\leq
\theta(G)$ for all graphs and $\theta(G)=\alpha(G)$ if $G$ is perfect.

Let $\tau(G)$ denote the size of a smallest vertex cover of $G$. For a
bipartite graph $G=(V,E)$, let $C$ be a minimum vertex cover of $G$ and $M$ a
maximum matching of $G$ that is formed from $C$ and a neighbor of each vertex
in $C$. If the end-vertices of $M$, $M_{c}$, yield the set $V$, then
$\theta_{c}(G)=\alpha(G)=|M|=\tau(G)$ and we are done. Otherwise,
$|M_{c}|<|V|$. Let $M_{u}$ be $V-M_{c}$.

\begin{proposition}
Let $G$ be a bipartite graph. Then $\theta_{c}(G) \leq\tau(G) + |M_{u}| =
\alpha(G)$.
\end{proposition}

\noindent\textbf{Proof}.\hspace{0.1in} Observe that $\alpha(G) = \tau(G) +
|M_{u}|$. Partition $V$ into sets such that each set contains the two
end-vertices from one edge in $M$; each vertex in $M_{u}$ is placed in a set
with a neighbor (which is a vertex in $M_{c}$). Note that each such set
induces a star. From this partitioning, we see that a neo-colonization exists
consisting only of stars -- and a star that is a $K_{2}$ has weight one and a
star that is a $K_{1, m}, m > 1$ has weight two. Therefore $\theta_{c}(G)
\leq\tau(G) + |M_{u}|$. \hfill$\blacksquare$

\medskip

As shown in \cite{KM2}, $\gamma_{\mathrm{m}}^{\infty}(T)=\theta_{c}(T)$ for
all trees $T$. There exist graphs with $\gamma(G)=\gamma_{\mathrm{m}}^{\infty
}(G)<\alpha(G)$, such as $C_{4}$ with a pendant vertex attached to one of its
vertices. Additional results comparing the vertex cover and eternal domination
numbers can be found in \cite{KM6}.

The following fact and its converse for $k=2$ (Proposition \ref{Prop_gamma2})
can be useful.

\begin{Fact}
\label{Fact_Necessary}A necessary condition for $\gamma^{\infty}(G)=k$, or
$\gamma_{\mathrm{m}}^{\infty}(G)=k$, is that every vertex of $G$ be contained
in a dominating set of size $k$.
\end{Fact}

If $k=1$, then this condition is also sufficient, and if $k\geq3$, then it is
not sufficient: let $T$ be the tree obtained by joining a new leaf to each
stem of $P_{3k-4}$. Then every vertex of $T$ is contained in some dominating
set of size $k$, but $\gamma^{\infty}(T)=\gamma_{\mathrm{m}}^{\infty}(T)>k$
(first attack one leaf, then attack another leaf at distance $3k-5$ from the
first leaf). For $k=2$, the condition is not sufficient for $\gamma^{\infty}$
(if $G=K_{m,n}$, $n\geq m\geq3$, then any pair of vertices from different
partite sets form a dominating set, but $\gamma^{\infty}(G)=n$). We show that
it is sufficient for~$\gamma_{\mathrm{m}}^{\infty}$.

\begin{proposition}
\label{Prop_gamma2}If every vertex of the graph $G\neq K_{n}$ is contained in
a dominating set of size $2$, then $\gamma_{\mathrm{m}}^{\infty}(G)=2$.
\end{proposition}

\noindent\textbf{Proof}.\hspace{0.1in} Suppose every vertex of $G$ is in a
dominating set of size two. Let $D=\{u,v\}$ be any dominating set and consider
any $x\in V-\{u,v\}$. We need to show that guards occupying $u$ and $v$ can
move to $x$ and to a vertex $y$ such that $\{x,y\}$ is a dominating set; that
is, $G$ has a dominating set $\{x,y\}$ such that $ux\in E(G)$ and $v\in N[y]$,
or $vx\in E(G)$ and $u\in N[y]$. Since $D$ dominates $x$, assume without loss
of generality that $vx\in E(G)$. By the hypothesis there exists a vertex $y$
such that $D^{\prime}=\{x,y\}$ is a dominating set. If $y\in N[u]$, we are
done. If $y\notin N[u]$, then $y\in N[v]$ because $D$ dominates $y$, and $u\in
N[x]$ because $D^{\prime}$ dominates $u$. But then $ux\in E(G)$ and $v\in
N[y]$, as required.\hfill$\blacksquare$

\section{m-Eternal Domination and Independence}

\label{Sec_Small_alpha}Clearly, if $\alpha(G)=1$ or $2$, then $\gamma
_{\mathrm{m}}^{\infty}(G)=\alpha(G)$. We next examine graphs with independence
number three, in which case $\gamma_{\mathrm{m}}^{\infty}(G)\in\{2,3\}$
(Theorem \ref{GHH1}). Classifying the graphs with $\alpha(G)=3$ and
$\gamma_{\mathrm{m}}^{\infty}(G)=2$, or equivalently $\alpha(G)=3$ and
$\gamma_{\mathrm{m}}^{\infty}(G)=3$, will make a valuable contribution to the
study of graphs with $\gamma_{\mathrm{m}}^{\infty}(G)=\alpha(G)$, but even
this apparently \textquotedblleft small\textquotedblright\ case may be
difficult as there is no known characterization of graphs with $\gamma=2$ and
$\alpha=3$.%
\begin{figure}[ptb]%
\centering
\includegraphics[
height=1.2073in,
width=2.0764in
]%
{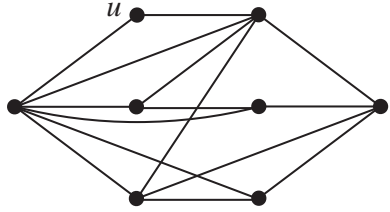}%
\caption{$\gamma_{\mathrm{m}}^{\infty}(G)=\alpha(G)=3$}%
\label{2-m-eternal}%
\end{figure}

The statement \textquotedblleft$\alpha(G)=3$ and any three independent
vertices of $G$ have a common neighbor\textquotedblright\ does not imply that
$\gamma_{\mathrm{m}}^{\infty}(G)=2$: for the graph $G$ in
Fig.~\ref{2-m-eternal}, $\alpha(G)=3$ and any three independent vertices of
$G$ have a common neighbor. However, the vertex $u$ is not in any dominating
set of size two. By Fact \ref{Fact_Necessary}, $\gamma_{\mathrm{m}}^{\infty
}(G)>2$, hence by Theorem \ref{GHH1}, $\gamma_{\mathrm{m}}^{\infty}(G)=3$.

We need to impose a stronger condition for the next result.

\begin{proposition}
\label{m-two-converse} Let $G=(V,E)$ be a graph with $\alpha(G)=3$. If $G$ has
a vertex $v$ that dominates all three vertices in all maximum independent
sets, then $\gamma_{\mathrm{m}}^{\infty}(G)=2$.
\end{proposition}

\noindent\textbf{Proof}.\hspace{0.1in} Since $\alpha(G)=3$, $\gamma
_{\mathrm{m}}^{\infty}(G)\geq2$. If $N[v]=V$ and $u\in V-\{v\}$ is arbitrary,
then $\{u,v\}$ is a domination set and the result follows from Proposition
\ref{Prop_gamma2}. Hence assume $X=\overline{N[v]}\neq\varnothing$. For any
distinct $x,x^{\prime}\in X$, $xx^{\prime}\in E$, otherwise $\{v,x,x^{\prime
}\}$ is an independent set not dominated by $v$. Thus $X$ is a clique. For any
$x\in X$ and any two distinct vertices $u,w\in\overline{N[x]}$, $uw\in E(G)$,
otherwise $\{x,u,w\}$ is an independent set not dominated by $v$; that is,
$\overline{N[x]}$ is a clique. Since $X$ is a clique, $\{v,x\}$ dominates $G$
for any $x\in X$. For any $u\in N(v)$, if $u$ is adjacent to all vertices in
$X$, then $\{u,v\}$ dominates $G$, and if $u$ is nonadjacent to some $x\in X$,
then the fact that $\overline{N[x]}$ is a clique implies that $\{u,x\}$
dominates $G$. Hence each vertex of $G$ is contained in a dominating set of
size two, and by Proposition \ref{Prop_gamma2}, $\gamma_{\mathrm{m}}^{\infty
}(G)=2$.\hfill$\blacksquare$

\bigskip

Note that $\gamma_{\mathrm{m}}^{\infty}(C_{6})=2$, $\alpha(C_{6})=3$, and no
maximum independent set is dominated by a single vertex. This example can be
generalized as follows to obtain a class of graphs $G$ such that
$\gamma_{\mathrm{m}}^{\infty}(G)=2$ and $\alpha(G)=3$. In $C_{6}=v_{0}%
,v_{1},...,v_{5},v_{0}$, replace each $v_{i}$ by a complete graph $H_{i}$ of
any order, and join each vertex of $H_{i}$, $i=0,...,5$, to each vertex of
$H_{i+1\ (\operatorname{mod}\ 6)}$ and to each vertex of
$H_{i-1\ (\operatorname{mod}\ 6)}$ to form the graph $H$. Note that
$\alpha(H)=3$ and, by Proposition \ref{Prop_gamma2}, $\gamma_{\mathrm{m}%
}^{\infty}(H)=2$ -- for any $u\in H_{i}$ and any $v\in
H_{i+3\ (\operatorname{mod}\ 3)}$, $\{u,v\}$ dominates $H$, $i=0,...,5$. Any
graph $G$ with $\alpha(G)=3$ that has $H$ as spanning subgraph also has
$\gamma_{\mathrm{m}}^{\infty}(G)=2$.

\section{Bipartite Graphs with $\gamma=\gamma^{\infty}$ or $\gamma
=\gamma_{\mathrm{m}}^{\infty}$}

\label{Sec2}In this section we consider bipartite graphs $G$ such that
$\gamma(G)=\gamma^{\infty}(G)$ or $\gamma(G)=\gamma_{\mathrm{m}}^{\infty}(G)$.
The former condition is more restrictive and this class of graphs is easy to
characterize. The second class is larger and more difficult to characterize,
and as a first step in this investigation we impose the further condition that
$\gamma(G)=2$. Recall that a graph is \textit{well-covered} if every maximal
independent set is maximum independent. For a matching $M$ in $G$, let $M(x)$
denote the vertex matched with $x$.

\begin{theorem}
\emph{\cite{rav}} \label{wc} A bipartite graph $G$ without isolated vertices
is well-covered if and only if $G$ has a perfect matching $M$ such that, for
every pair $(x,M(x))$, the subgraph induced by $N(x)\cup N(M(x))$ is complete bipartite.
\end{theorem}

\begin{proposition}
\label{PropBipartite}Let $G$ be a bipartite graph without isolated vertices.
Then $\gamma(G)=\gamma^{\infty}(G)$ if and only if $\gamma(G)=n/2$.
\end{proposition}

\noindent\textbf{Proof}.\hspace{0.1in} If $\gamma(G)=n/2$, then $G$ is
well-covered. By Theorem \ref{wc}, $G$ has a perfect matching. Since $G$ is
bipartite, $\theta(G)=n/2$, which implies $\gamma^{\infty}(G)=n/2$. On the
other hand, if $\gamma(G)=\gamma^{\infty}(G)$, then, by Fact
\ref{FactED_Bound}, $\gamma^{\infty}(G)=\alpha(G)$. Since $\alpha(G)\geq n/2$
for any bipartite graph, the result follows. \hfill$\blacksquare$

\bigskip

Note that $\gamma(G)=n/2$ if and only if each component of $G$ is a 4-cycle or
the corona of a connected graph $H$ with $K_{1}$, c.f. \cite{HHS}. We
strengthen Proposition \ref{PropBipartite} to triangle-free graphs in
Corollary \ref{CorK3_free}.

If $\gamma(G)=1$, then $\gamma_{\mathrm{m}}^{\infty}(G)=1$ if $G$ is complete,
and $\gamma_{\mathrm{m}}^{\infty}(G)=2$ otherwise. Now we turn to describing
the bipartite graphs with $\gamma_{\mathrm{m}}^{\infty}=\gamma=2$. Let
$\mathfrak{C}$ be the class of all graphs obtained from $K_{m,m}$, $m\geq2$,
by deleting a matching $M$ of size $k$, where $0\leq k\leq m$, or from
$K_{m,n}$, $n>m\geq2$, by deleting a matching $M$ of size $\ell$, $0\leq
\ell\leq m-1$. For example, $\mathfrak{C}$ contains the graphs $2K_{2}%
,\ P_{4},\ C_{6},\ K_{m,n},\ K_{2,3}-e$. If $G\in\mathfrak{C}$ and $v$ is a
vertex of $G$ incident with an edge of the removed matching $M$, then $v$ is a
\emph{depleted vertex}, otherwise $v$ is a \emph{full vertex}. Note that each
$G\in\mathfrak{C}$ that has a full vertex, has a full vertex in each of its
partite sets.

\begin{theorem}
If $G$ is bipartite, then $\gamma(G)=\gamma_{\mathrm{m}}^{\infty}(G)=2$ if and
only if $G\in\mathfrak{C}$.
\end{theorem}

\noindent\textbf{Proof}.\hspace{0.1in} Let $G$ have partite sets $A$ and $B$.
Suppose $G\in\mathfrak{C}$. Then $\gamma_{\mathrm{m}}^{\infty}(G)\geq
\gamma(G)\geq2$. If $x\in A$ is full, then there exists $y\in B$ that is full,
and $\{x,y\}$ dominates $G$. If $x\in A$ is depleted, let $y\in B$ be the
vertex such that $xy$ belongs to the deleted matching. Then $\{x,y\}$
dominates $G$. Hence each vertex of $A$, and similarly each vertex of $B$
belongs to a dominating set of size two. By Proposition \ref{Prop_gamma2},
$\gamma_{\mathrm{m}}^{\infty}(G)=\gamma(G)=2$.

Conversely, suppose $\gamma(G)=\gamma_{\mathrm{m}}^{\infty}(G)=2$. Then $G$
does not have a universal vertex, so $|A|,|B|\geq2$. Assume without loss of
generality that $2\leq m=|A|\leq n=|B|$.

Suppose $\deg v\leq m-2$ for some $v\in B$; say $v$ is nonadjacent to
$u,u^{\prime}\in A$. By Fact \ref{Fact_Necessary} there is a configuration of
guards such that $u$ is occupied. Since $u^{\prime}$ is protected, the other
guard occupies $u^{\prime}$ or some vertex $w\in B-\{v\}$. But in either case
$v$ is unprotected, contradicting $\gamma_{\mathrm{m}}^{\infty}(G)=2$. Hence
$\deg v\geq m-1$ for each $v\in B$. Similarly, $\deg u\geq n-1$ for each $u\in
A$. Therefore $G=K_{m,n}$ or $G$ is obtained from $K_{m,n}$ by deleting edges
of a matching.

Now suppose $m<n$ and $\deg u=n-1$ for each $u\in A$. Since $m<n$ there exists
$v\in B$ such that $\deg v=m$. Let $v$ be occupied. Since $|B-\{v\}|\geq2$,
the other guard occupies a vertex $u\in A$. Now $v$ is adjacent to $u$, and
$\deg u=n-1$; hence there exists $w\in B-\{v\}$ such that $uw\notin E(G)$. But
then $w$ is not protected, a contradiction as above. We deduce that $\deg u=n$
for at least one vertex $u\in A$. Therefore $G\in\mathfrak{C}$ as
required.\hfill$\blacksquare$

\bigskip

It turns out that the class of triangle-free graphs with $\gamma_{\mathrm{m}%
}^{\infty}=\gamma=2$ is almost the same as the class of bipartite graphs with
this property.

\begin{corollary}
\label{CorTriangle-free}A triangle-free graph $G$ satisfies $\gamma
(G)=\gamma_{\mathrm{m}}^{\infty}(G)=2$ if and only if $G=C_{5}$ or
$G\in\mathfrak{C}$.
\end{corollary}

\noindent\textbf{Proof}.\hspace{0.1in} Suppose $G\ncong C_{5}$ is a
non-bipartite triangle-free graph such that $\gamma(G)=\gamma_{\mathrm{m}%
}^{\infty}(G)=2$. Then $G$ has a shortest odd cycle $H\cong C_{2n+1}$, where
$n\geq2$. Since the component of $G$ containing $H$ is not complete and
$\gamma(G)=2$, $G$ is connected. We obtain a contradiction by proving by
induction on $n$ that $H\ncong C_{2n+1}$ for all $n\geq2$.

Suppose first that $H\cong C_{5}$; say $H$ is the cycle $v_{0},v_{1}%
,...,v_{4},v_{0}$. Since $H$ is triangle-free, $H$ is a chordless $5$-cycle.
Since $G\ncong C_{5}$, there exists a vertex $x\in V(G)-V(H)$ that is adjacent
to a vertex of $H$; say $xv_{0}\in E(G)$. By Fact \ref{Fact_Necessary} there
exists a vertex $y$ such that $\{x,y\}$ is a dominating set of $G$.

Suppose $x$ is not adjacent to any other vertex of $H$. Then $y$ dominates
$\{v_{1},...,v_{4}\}$. Since $G[\{v_{1},...,v_{4}\}]\cong P_{4}$ and no vertex
of $G-H$ dominates more than two of $v_{1},...,v_{4}$, this is impossible.
Hence $x$ is adjacent to $v_{i}$ for some $i=1,...,4$. Since $G$ is
triangle-free, we may assume without loss of generality that $xv_{2}\in E(G)$
and $xv_{i}\notin E(G)$ for $i=1,3,4$. Then $y$ dominates $\{v_{1},v_{3}%
,v_{4}\}$. But $G$ is triangle-free, so neither $v_{3}$ nor $v_{4}$ dominates
$v_{1}$, and no other vertex of $G$ dominates both $v_{3}$ and $v_{4}$. We
deduce that $H\ncong C_{5}$.

Now suppose that for some $k\geq3$, $H\ncong C_{2r+1}$ for all $r=2,...,k-1$
and suppose $H\cong C_{2k+1}$. Say $H$ is the cycle $v_{0},v_{1}%
,...,v_{2k},v_{0}$. Since $\gamma(C_{2k+1})>2$, $G\ncong H$. If $H$ has a
chord, then $G$ has an odd cycle $C_{2r+1}$ for $r<k$, which is not the case.
Hence there is a vertex $x\in V(G)-V(H)$ such that $x$ is adjacent to a vertex
of $H$, say to $v_{0}$. As before, there is a vertex $y$ such that $\{x,y\}$
is a dominating set of $G$. If $x$ is not adjacent to any other vertex of $H$,
we obtain a contradiction as in the case where $H=C_{5}$. On the other hand,
if $x$ is adjacent to some $v_{j}$, $j\in\{1,...,2k\}-\{2,2k-1\}$, then $G$
also has an odd cycle $C_{2r+1}$ for $r<k$. Hence assume $xv_{2}\in E(G)$.
Then $x$ is not adjacent to $v_{2k-1}$, hence $y$ dominates all of
$v_{1},v_{3},v_{4},...,v_{2k}$. As in the case where $H=C_{5}$, this is impossible.

By induction, $H\ncong C_{2n+1}$ for all $n\geq2$. Therefore $C_{5}$ is the
only non-bipartite triangle-free graph $G$ such that $\gamma(G)=\gamma
_{\mathrm{m}}^{\infty}(G)=2$. \hfill$\blacksquare$

\section{Trees with $\gamma_{\mathrm{m}}^{\infty}=\theta$}

\label{SecTrees}In this section we prove our main result -- a characterization
of the class of trees $T$ for which $\gamma_{\mathrm{m}}^{\infty}%
(T)=\theta(T)$. We begin by stating two reductions on trees from \cite{KM2}.

\smallskip\noindent\textbf{R1}: Let $x$ be a stem of $T$ adjacent to $\ell
\geq2$ leaves and to exactly one vertex of degree at least two. Delete all
leaves adjacent to $x$.

\smallskip\noindent\textbf{R2}: Let $x$ be a stem of degree two in $T$ such
that $x$ is adjacent to exactly one leaf, $y$. Delete both $x$ and $y$.

\begin{lemma}
\emph{\cite{KM2}} \label{rlemma} If $T^{\prime}$ is the result of applying
reduction R1 or R2 to the tree $T$, then $T^{\prime}$ is a tree and
$\gamma_{\mathrm{m}}^{\infty}(T)=1+\gamma_{\mathrm{m}}^{\infty}(T^{\prime})$.
\end{lemma}

It is shown in \cite{KM2} that one can repeatedly apply these reductions,
reducing $T$ to a star $K_{1,r}$, $r\geq1$, in such a way as to compute
$\theta_{c}(T)=\gamma_{\mathrm{m}}^{\infty}(T)$. The characterization of trees
with equal clique covering and $\mathrm{m}$-eternal domination numbers follows.

\begin{theorem}
\label{ThmTheta}Let $T$ be a tree with at least two vertices. Then
$\gamma_{\mathrm{m}}^{\infty}(T)=\theta(T)$ if and only if the reduction R2
can be applied repeatedly to $T$ to obtain a star $K_{1,r}$, $r\in\{1,2\}$.
\end{theorem}

\noindent\textbf{Proof}.\hspace{0.1in} Suppose first that $T=K_{1,r}$,
$r\geq1$. Then either $T=K_{2}$ and $\gamma_{\mathrm{m}}^{\infty}%
(T)=\theta(T)=1$, or $r\geq2$, $\gamma_{\mathrm{m}}^{\infty}(T)=2$ and
$\theta(T)=r$, hence $\gamma_{\mathrm{m}}^{\infty}(T)=\theta(T)=2$ if and only
if $r=2$ and thus $T=K_{1,2}$. Hence the theorem holds for stars. Assume the
theorem holds for all trees of order less than $n$, where $n\geq4$, and let
$T$ be a tree of order $n$. We may assume that $T$ is not a star.

First assume that $T$ can be reduced to $K_{2}$ or $K_{1,2}$ by repeatedly
applying R2. Since $T$ is not a star, $T$ has a stem $x$ of degree two that is
adjacent to exactly one leaf, say $y$, such that $T^{\prime}=T-\{x,y\}$ is
either $K_{2}$, $K_{1,2}$ or can be reduced to one of these trees by
repeatedly applying R2. By the induction hypothesis, $\gamma_{\mathrm{m}%
}^{\infty}(T^{\prime})=\theta(T^{\prime})$. By Lemma \ref{rlemma},
$\gamma_{\mathrm{m}}^{\infty}(T)=1+\gamma_{\mathrm{m}}^{\infty}(T^{\prime})$,
and obviously $\theta(T)=\theta(T^{\prime})+1$, so that $\gamma_{\mathrm{m}%
}^{\infty}(T)=\theta(T)$.

Conversely, assume $T$ cannot be reduced to $K_{2}$ or $K_{1,2}$ by repeatedly
applying R2. Apply R2 to $T$ repeatedly until a tree $T^{\prime}\notin%
\{K_{2},K_{1,2}\}$ is obtained to which R2 cannot be applied; say R2 is
applied $k$ times to obtain $T^{\prime}$. By Lemma \ref{rlemma} applied $k$
times, $\gamma_{\mathrm{m}}^{\infty}(T^{\prime})=\gamma_{\mathrm{m}}^{\infty
}(T)-k$. Similarly, each application of R2 reduces the clique partition number
by 1, thus $\theta(T^{\prime})=\theta(T)-k$. Therefore, if we can show that
$\gamma_{\mathrm{m}}^{\infty}(T^{\prime})<\theta(T^{\prime})$, it will follow
that $\gamma_{\mathrm{m}}^{\infty}(T)<\theta(T)$ and the proof will be
complete. The remainder of the proof shows that $\theta_{c}(T^{\prime}%
)<\theta(T^{\prime})$.

If $T^{\prime}$ is a star, then $T^{\prime}=K_{1,r}$, $r\geq3$, and
$\theta_{c}(T^{\prime})=\gamma_{\mathrm{m}}^{\infty}(T^{\prime})=2<r=\theta
(T^{\prime})$. Hence assume $T^{\prime}$ is not a star. Since R2 cannot be
performed on $T^{\prime}$, each stem of $T^{\prime}$ is a branch vertex and
$T^{\prime}$ has at least two branch vertices, hence at least two end-branch
vertices. Moreover, each end-branch vertex $v$ is adjacent to $\deg v-1$
leaves and one non-leaf vertex of $T^{\prime}$. Note that each clique
partition of $T^{\prime}$ is a neo-colonization. Consider a minimum clique
partition $\Theta=\{U_{0},...,U_{\theta-1}\}$ of $T^{\prime}$ (thus each
$U_{i}$ induces a $K_{1}$ or a $K_{2}$). We show that there exists a
neo-colonization $\mathcal{P}$ of $T^{\prime}$ with $w(\mathcal{P})<w(\Theta
)$. The result $\theta_{c}(T^{\prime})<\theta(T^{\prime})$ then follows.

\begin{figure}[ptb]%
\centering
\includegraphics[
height=1.625in,
width=2.7622in
]%
{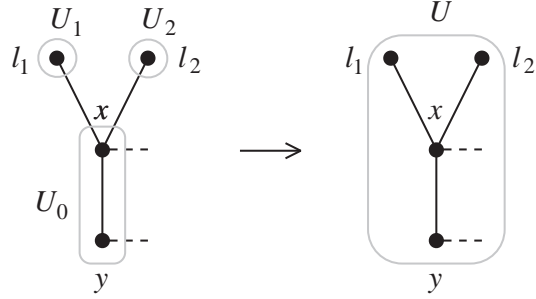}%
\caption{$w(U_{0})+w(U_{1})+w(U_{2})=3$ and $w(U)=\gamma_{c}(K_{1,3})+1=2$}%
\label{Fig5.2a}%
\end{figure}

Suppose $T^{\prime}$ has a stem $x$ adjacent to leaves $\ell_{1}$ and
$\ell_{2}$ such that $\{\ell_{i}\}$ is a part of $\Theta$ for $i=1,2$; without
loss of generality say $U_{i}=\{\ell_{i}\}$, $i=1,2$. See Fig.~\ref{Fig5.2a}.
Since $\Theta$ is a minimum clique cover, there exists $y\in N(x)-\{\ell
_{1},\ell_{2}\}$ such that $\{x,y\}$ is a part of $\Theta$; say $U_{0}%
=\{x,y\}$. Then $w(U_{i})=1$, $i=0,1,2$. Let $U=\bigcup_{i=0}^{2}U_{i}$ and
note that $T^{\prime}[U]=K_{1,3}$. Let $\mathcal{P}$ be the neo-colonization
of $T^{\prime}$ defined by $\mathcal{P}=(\Theta-\{U_{0},U_{1},U_{2}%
\})\cup\{U\}$ and note that $w(U)=\gamma_{c}(T^{\prime}[U])+1=2$. Then
$w(\mathcal{P})=w(\Theta)-3+2=w(\Theta)-1<\theta(T^{\prime})$ and we are done.
Hence we may assume that each stem of $T^{\prime}$ is adjacent to at most one
leaf $\ell$ such that $U_{i}=\{\ell\}$ for some $i$. In particular, each
end-branch vertex $x$ has degree three and is adjacent to leaves $x_{1},x_{2}$
such that (say) $\{x_{1}\}$ and $\{x,x_{2}\}$ are parts of $\Theta$.%

\begin{figure}[ptb]%
\centering
\includegraphics[
height=3.0761in,
width=5.9369in
]%
{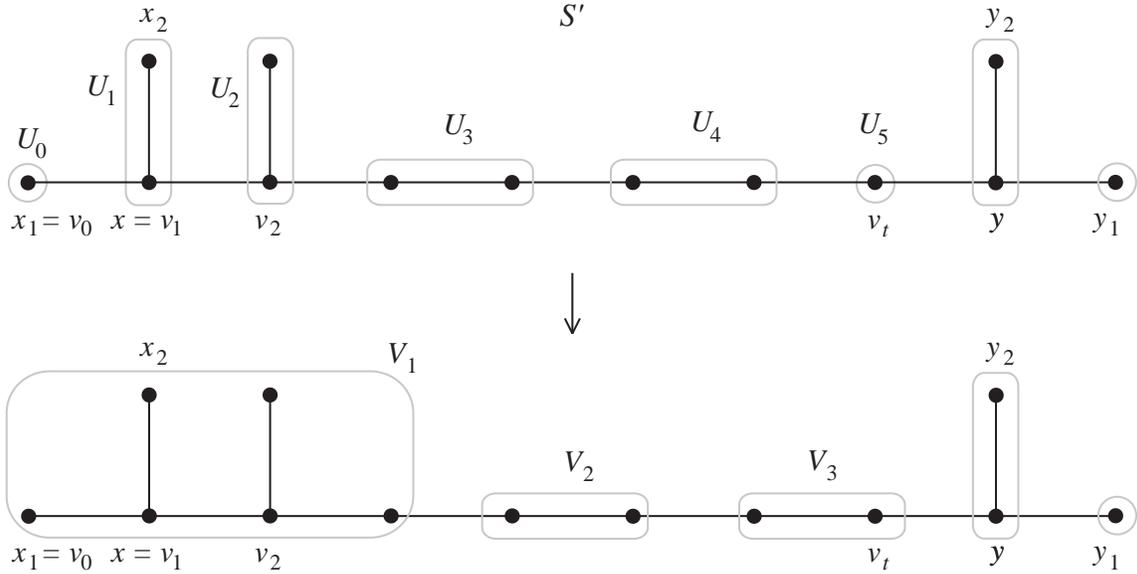}%
\caption{$\sum_{i=0}^{5}w(U_{i})=6$ and $\sum_{i=1}^{3}w(V_{i})=5$}%
\label{Fig5.2b}%
\end{figure}

Let $x$ and $y$ be two end-branch vertices of $T^{\prime}$, with $x_{1}$ and
$x_{2}$ as above, and let $y_{1},y_{2}$ be the leaves adjacent to $y$ such
that $\{y_{1}\}$ and $\{y,y_{2}\}$ are parts of $\Theta$. Let $Q^{\prime
}:x_{1}=v_{0},...,v_{t^{\prime}}=y_{1}$ be the $x_{1}$-$y_{1}$ path in
$T^{\prime}$. (Thus $v_{1}=x$ and $v_{t^{\prime}-1}=y$.) With respect to
$Q^{\prime}$, we consider three types of parts $U_{i}$ of $\Theta$: a $K_{1}%
$-part $\{u\}$, where $u\in V(Q^{\prime})$, a part $\{u,u^{\prime}\}$, where
$u,u^{\prime}\in V(Q^{\prime})$, which we refer to as a $K_{2}$-part, and a
part $\{u,u^{\prime}\}$, where $\{u,u^{\prime}\}\cap V(Q^{\prime})=\{u\}$,
which we refer to as a $P_{2}$-part. Since $\{v_{t^{\prime}}\}$ is a $K_{1}%
$-part on $Q^{\prime}$, there exists a smallest integer $t$, $1\leq t\leq
t^{\prime}$, such that $\{v_{t}\}$ is a $K_{1}$-part on $Q^{\prime}$. Let
$Q:v_{0},...,v_{t}$ be the $v_{0}$-$v_{t}$ subpath of $Q^{\prime}$. Note that
$\{x,x_{2}\}$ is a $P_{2}$-part. Therefore the parts $\Omega=\{U_{i}:U_{i}\cap
V(Q)\neq\varnothing\}$ of $\Theta$ form a sequence that consists of a $K_{1}%
$-part $\{v_{0}\}=\{x_{1}\}$, followed by a number of $P_{2}$ parts, followed
(possibly) by a number of $K_{2}$-parts, then $P_{2}$-parts, and so on,
finally ending in the $K_{1}$-part $\{v_{t}\}$. We can therefore define a
sequence of positive integers $s_{1},s_{2},...,s_{k}$ such that the part
$\{v_{0}\}$ is followed by $s_{1}$ $P_{2}$-parts, the last of which is
followed by $s_{2}$ $K_{2}$-parts, then $s_{3}$ $P_{2}$-parts, and so on,
until the final $s_{k}$ $K_{2}$- or $P_{2}$-parts are followed by $\{v_{t}\}$.
See the top graph in Fig.~\ref{Fig5.2b}. Let $\omega=w(\Omega)$. Since each
part of $\Theta$ is assigned a weight of one when $\Theta$ is considered as a
neo-colonization,
\begin{equation}
\omega=w(\Omega)=2+\sum_{i=1}^{k}s_{i}. \label{eqOmega}%
\end{equation}
We may assume that the parts of $\Theta$ that belong to $\Omega$ are labeled
$U_{0}=\{v_{0}\},\ U_{1}=\{v_{1},x_{2}\},\ ...,\ U_{s_{1}},\ U_{s_{1}%
+1},\ ...,\ U_{s_{1}+s_{2}},\ ...,\ U_{\omega}=\{v_{t}\}$, in order of their
occurrence on $Q$. Thus $U_{1},...,U_{s_{1}}$ are $P_{2}$-parts, $U_{s_{1}%
+1},...,U_{s_{1}+s_{2}}$ are $K_{2}$ parts, and so on. Let $S^{\prime}$ be the
subgraph of $T^{\prime}$ induced by $\bigcup_{i=0}^{\omega}U_{i}$. Since
$\Theta$ is a clique cover of $T$ and each vertex of $Q$ is contained in a set
$U_{i}$, $i=0,...,\omega$, $S^{\prime}$ is a tree. We define a
neo-colonization $\mathcal{P}^{\prime}=\{V_{1},...,V_{r}\}$ of $S^{\prime}$ as follows.%

As illustrated in Fig.~\ref{Fig5.2b}, we combine each subsequence of
consecutive $P_{2}$-parts with the last vertex of $Q$ preceding and the first
vertex of $Q$ following this subsequence into one part. We also combine the
second vertex of each $K_{2}$-part with the first vertex of the next $K_{2}%
$-part to form new $K_{2}$-parts, ending with $v_{t}$ belonging to either a
$K_{2}$-part or a part containing $P_{2}$-parts. In order to calculate the
weight of $\mathcal{P}^{\prime}$, we describe the process more formally.

\begin{itemize}
\item Let $V_{1}$ consist of $\bigcup_{j=0}^{s_{1}}U_{j}$ together with the
first vertex of $Q$ that belongs to $U_{s_{1}+1}$. Then $S^{\prime}[V_{1}]$ is
connected and $\gamma_{c}(S^{\prime}[V_{1}])=s_{1}$.

\item For $i=2,...,s_{2}$, let $V_{i}$ consist of the second vertex of $Q$
that belongs to $U_{s_{1}+i-1}$ and the first vertex of $Q$ that belongs to
$U_{s_{1}+i}$; each such $V_{i}$ is a $K_{2}$-part.

\item If $v_{t}$ has not been reached above, let $V_{s_{2}+1}$ consist of
$\bigcup_{j=s_{2}+1}^{s_{3}}U_{j}$ together with the last vertex of $Q$ that
belongs to $U_{s_{2}}$ and the first vertex of $Q$ that belongs to
$U_{s_{3}+1}$. Then $S^{\prime}[V_{_{s_{2}+1}}]$ is connected and $\gamma
_{c}(S^{\prime}[V_{_{s_{2}+1}}])=s_{3}$.

\item Continue by splitting and recombining the next $K_{2}$-parts, if necessary.

\item Finally, $V_{r}$ either consists of $v_{t}$ and the last vertex of
$U_{\omega-1}$, if $U_{\omega-1}$ is a $K_{2}$-part, or of the union of the
last $s_{k}$ consecutive $P_{2}$-parts of $\Omega$ on $Q$, together with
$v_{t}$ and the last vertex of $Q$ that belongs to $U_{s_{1}+\cdots
+s_{\omega-2}}$, otherwise.
\end{itemize}

The sets $V_{i}$ are mutually disjoint, each $S^{\prime}[V_{i}]$ is connected
and $\bigcup_{i=1}^{r}V_{i}=V(S^{\prime})$. Hence $\mathcal{P}^{\prime}$ is a
neo-colonization of $S^{\prime}$. The weight $w(\mathcal{P}^{\prime})$ is
calculated as follows. If $V_{i}$ contains $s_{j}$ $P_{2}$-parts, then
$w(V_{i})=\gamma_{c}(S^{\prime}[V_{i}])+1=s_{j}+1$. Each such $V_{i}$, $i\neq
r$, is followed by $s_{j+1}-1$ $K_{2}$-parts of $\mathcal{P}^{\prime}$.
Therefore, if $V_{r}$ is a $K_{2}$-part of $\mathcal{P}^{\prime}$, then $k$ is
even, $w(V_{r})=1$ and
\[
w(\mathcal{P}^{\prime})=(s_{1}+1)+(s_{2}-1)+\cdots+(s_{k-1}+1)+(s_{k}%
-1)+1=1+\sum_{i=1}^{k}s_{i}=w(\Omega)-1\ \ \text{(by (\ref{eqOmega}))},
\]
and if $V_{r}$ contains $P_{2}$-parts of $\Omega$, then $k$ is odd and, again
using (\ref{eqOmega}),%
\[
w(\mathcal{P}^{\prime})=(s_{1}+1)+(s_{2}-1)+\cdots+(s_{k-1}-1)+(s_{k}%
+1)=1+\sum_{i=1}^{k}s_{i}=w(\Omega)-1.
\]
Let $\mathcal{P}=\mathcal{P}^{\prime}\cup\{U_{i}\in\Theta:U_{i}\cap
V(S^{\prime})=\varnothing\}$. Then $\mathcal{P}$ is a neo-colonization of
$T^{\prime}$ and%
\[
w(\mathcal{P})=w(\mathcal{P}^{\prime})+w(\Theta-\Omega)\leq w(\Omega
)-1+w(\Theta)-w(\Omega)<w(\Theta).
\]
Therefore $\theta_{c}(T^{\prime})<\theta(T^{\prime})$, hence $\gamma
_{\mathrm{m}}^{\infty}(T)=\theta_{c}(T^{\prime})<\theta(T^{\prime})$%
.\hfill$\blacksquare$

\bigskip

The next result follows immediately from Theorem \ref{ThmTheta}.

\begin{corollary}
If $T$ is a tree with at least two vertices, then $\gamma_{\mathrm{m}}%
^{\infty}(T)=\theta(T)$ if and only if $T$ can be obtained from $K_{2}$ or
$P_{3}$ by successively adding a new $K_{2}$, joining one of its leaves to any
vertex of the previously constructed tree.
\end{corollary}

\section{Clique Covering Numbers of Graphs with $\gamma=\gamma^{\infty}$}

\label{SecTheta}There are many graphs with $\gamma(G)=\theta(G)$, including
$C_{4}$, and two $K_{n}$'s connected by one edge, though the two parameters
may also differ by any arbitrary amount, for example in $K_{1,m}$. There does
not exist a meaningful characterization of the graphs $G$ with $\gamma
(G)=\theta(G)$ and this complicates the issue of characterizing graphs with
$\gamma^{\infty}(G)=\theta(G)$. The results of this section are motivated by
an error discovered in \cite{KM2}, where it was claimed that if $\gamma
(G)=\gamma^{\infty}(G)$, then $\gamma(G)=\theta(G)$. The proof given in
\cite{KM2} is incorrect, as the initial set of cliques consisting of the
vertices in dominating set $D$ and their private neighbors cannot, in fact, be
extended to other vertices of $G$. We determine two classes of graphs $G$ such
that $\gamma(G)=\gamma^{\infty}(G)=\theta(G)$. The following fact was proved
in \cite{KM2}, and will be needed below.

\begin{Fact}
\label{FactEDS}Let $D$ be an EDS of a graph $G$. For each $v\in D$,
$G[\{v\}\cup\operatorname{epn}(v,D)]$ is a clique, and if $v\in D$ defends
$u\in V(G)-D$, then $G[\{u,v\}\cup\operatorname{epn}(v,D)]$ is a clique.
\end{Fact}

As shown in \cite{BC}, every graph without isolated vertices has a minimum
dominating set $D$ such that $\operatorname{epn}(v,D)\neq\varnothing$ for each
$v\in D$. A similar result does not hold for minimum EDS's -- consider $P_{3}%
$, for example. We now prove a corresponding result under restricted
conditions. If $D$ is an EDS of a graph $G$, and $w\in V(G)-D$ is adjacent to
more than one vertex in $D$, we say that $w$ is a \emph{shared vertex}.

\begin{lemma}
\label{Lem_epn}If $G$ is a graph without isolated vertices such that
$\gamma(G)=\gamma^{\infty}(G)$ and $\Delta(G)\leq3$, then $G$ has a minimum
EDS $D$ such that $\operatorname{epn}(v,D)\neq\varnothing$ for each $v\in D$.
\end{lemma}

\noindent\textbf{Proof}.\hspace{0.1in} Let $D$ be a minimum EDS of $G$ that
maximizes the number of edges in $G[D]$. We first show that
\begin{equation}
\text{if}\ u\in D\ \text{and}\ \operatorname{epn}(u,D)=\varnothing
,\ \text{then}\ u\ \text{does\ not\ defend\ any\ vertex\ of}\ G-D.
\label{eqLem_pn1}%
\end{equation}
Suppose $u\in D$ and $\operatorname{epn}(u,D)=\varnothing$. Since
$\gamma(G)=\gamma^{\infty}(G)$, $D$ is a minimum dominating set, hence $u$ is
isolated in $G[D]$ (because $\operatorname{pn}(u,D)\neq\varnothing$). Suppose,
to the contrary, that $u$ defends $w\in V(G)-D$. Then $D^{\prime
}=(D-\{u\})\cup\{w\}$ is an EDS. Moreover, $w$ is adjacent to a vertex in
$D^{\prime}$, so that $G[D^{\prime}]$ has more edges than $G[D]$, contrary to
the choice of $D$.

Now we show that%
\begin{equation}
\text{each}\ w\in V(G)-D\ \text{is adjacent to at most two vertices in}\ D.
\label{eqLem_pn2}%
\end{equation}
Suppose $w\in V(G)-D$ is adjacent to more than two vertices in $D$. Since
$\Delta(G)\leq3$, $w$ is adjacent to exactly three vertices $v_{1},v_{2}%
,v_{3}\in D$ and nonadjacent to all external private neighbors of $v_{i}$,
$i=1,2,3$. But $D$ is an EDS, and some $v\in\{v_{1},v_{2},v_{3}\}$ defends
$w$. By Fact \ref{FactEDS}, $\operatorname{epn}(v,D)=\varnothing$. This
contradicts (\ref{eqLem_pn1}).

We also show that
\begin{align}
\text{if }u,v  &  \in D,\text{ }\operatorname{epn}(u,D)=\varnothing
,\ \text{and }u\text{ and }v\text{ have a shared neighbor in }G-D,\nonumber\\
&  \text{then they have exactly two shared neighbors in }G-D.
\label{eqLem_pn3}%
\end{align}
Suppose $N(u)\cap N(v)\cap(V-D)=\{w\}$. By (\ref{eqLem_pn2}), $N(w)\cap
D=\{u,v\}$. Since $D$ is an EDS and $u$ does not defend $w$ by
(\ref{eqLem_pn1}), $v$ defends $w$, $\operatorname{epn}(v,D)\neq\varnothing$,
and $w$ is adjacent to each vertex in $\operatorname{epn}(v,D)$ (Fact
\ref{FactEDS}). But then $(D-\{u,v\})\cup\{w\}$ dominates $G$, a contradiction
because $D$ is a minimum dominating set. On the other hand, suppose $N(u)\cap
N(v)\cap(V-D)=\{w_{1},w_{2},w_{3}\}$. Then $N(u)\cap N(v)=\{w_{1},w_{2}%
,w_{3}\}$ because $\Delta(G)\leq3$, hence $\operatorname{epn}%
(v,D)=\operatorname{epn}(u,D)=\varnothing$, and by (\ref{eqLem_pn1}), neither
$u$ nor $v$ defends $w_{i}$, $i=1,2,3$. But by (\ref{eqLem_pn2}),
$N(w_{i})\cap D=\{u,v\}$ and so no vertex in $D$ defends $w_{i}$, a contradiction.

Now consider $u\in D$ such that $\operatorname{epn}(u,D)=\varnothing$. As in
the proof of (\ref{eqLem_pn1}), $u$ is isolated in $G[D]$. Since
$\delta(G)\geq1$, $u$ has at least one neighbor in $G-D$. By (\ref{eqLem_pn3})
there exists $v\in D$ such that $N(u)\cap N(v)\cap(V-D)=\{w_{1},w_{2}\}$, say.
As in the proof of (\ref{eqLem_pn3}), $v$ defends $w_{1}$ and $w_{2}$, and
$\operatorname{epn}(v,D)\neq\varnothing$. Now $v$ is adjacent to three
vertices of $G-D$, hence $v$ is isolated in $G[D]$. Since $v$ defends $w_{1}$,
$D^{\prime\prime}=(D-\{v\})\cup\{w_{1}\}$ is an EDS. However, $w_{1}$ is
adjacent to $u$ in $D^{\prime\prime}$, which implies that $G[D^{\prime\prime
}]$ has more edges than $G[D]$, a contradiction.\hfill$\blacksquare$

\bigskip

We use Fact \ref{FactEDS} and Lemma \ref{Lem_epn} to prove the main result of
this section.

\begin{theorem}
\label{main} Let $G$ be a graph with $\gamma(G)=\gamma^{\infty}(G)$ and
$\Delta(G)\leq3$. Then $\gamma^{\infty}(G)=\theta(G)$.
\end{theorem}

\noindent\textbf{Proof}.\hspace{0.1in} We may assume without loss of
generality that $G$ has no isolated vertices. Let $D$ be a minimum EDS of $G$
such that $\operatorname{epn}(v,D)\neq\varnothing$ for each $v\in D$; such an
EDS exists by Lemma \ref{Lem_epn}. If $\gamma(G)=\gamma^{\infty}(G)=1$, then
$G$ is complete and the statement holds. Hence we assume $\gamma(G)>1$.

If each vertex of $G-D$ is an external private neighbor of a vertex in $D$,
then, by Fact \ref{FactEDS}, $\{\{x\}\cup\operatorname{epn}(x,D):x\in D\}$ is
a clique cover of $G$ and the result follows. Hence assume some vertex of
$G-D$ is a shared vertex. For each $x\in D$, let $S_{x}$ denote the set of
shared vertices defended by $x$. If $|S_{x}|\leq1$ for each $x\in D$, then
$R_{x}=\{x\}\cup S_{x}\cup\operatorname{epn}(x,D)$ forms a clique (Fact
\ref{FactEDS}) and $\{R_{x}:x\in D\}$ is a clique partition of $G$ into
$\gamma(G)$ parts.

Therefore we assume that $w,w^{\prime}\in S_{u}$ for some $u\in D$. Say $w$
and $w^{\prime}$ are also adjacent to $v$ and $v^{\prime}$, respectively,
where possibly $v=v^{\prime}$. Let $y\in\operatorname{epn}(v^{\prime},D)$ and
$z\in\operatorname{epn}(u,D)$. By Fact \ref{FactEDS}, $w$ and $w^{\prime}$ are
adjacent to $z$. Since $\Delta(G)\leq3$, $N(w)=\{u,v,z\}$ and $N(w^{\prime
})=\{u,v^{\prime},z\}$; note that $w,w^{\prime}$ are not adjacent to each
other or to $y$. Since $u$ defends $w$, $D^{\prime}=(D-\{u\})\cup\{w\}$ is an
EDS, and $\{w^{\prime},y\}\subseteq\operatorname{epn}(v^{\prime},D^{\prime})$.
Since $w^{\prime}$ is not adjacent to $y$, this contradicts Fact
\ref{FactEDS}.\hfill$\blacksquare$

\begin{corollary}
\label{CorK3_free}Let $G$ be a triangle-free graph such that $1\leq
\delta(G)\leq\Delta(G)\leq3$. Then $\gamma(G)=\gamma^{\infty}(G)$ if and only
if $\gamma(G)=n/2$.
\end{corollary}

\noindent\textbf{Proof}.\hspace{0.1in} Since $G$ has no isolated vertices,
$\gamma(G)\leq n/2$. Suppose $\gamma(G)=\gamma^{\infty}(G)$. By Theorem
\ref{main}, $\theta(G)=\gamma(G)$, and since $G$ is triangle-free,
$\theta(G)\geq n/2$. Conversely, suppose $\gamma(G)=n/2$ and let $D$ be a
minimum dominating set such that $\operatorname{epn}(v,D)\neq\varnothing$ for
each $v\in D$. Then $|\operatorname{epn}(v,D)|=1$ for each $v\in D$; say
$\operatorname{epn}(v,D)=\{v^{\prime}\}$. Then $\mathcal{P}=\{\{v,v^{\prime
}\}:v\in D\}$ is a clique partition of $G$. Since $G$ is triangle-free,
$\mathcal{P}$ is a minimum clique partition and so $\theta(G)=n/2$
.\hfill$\blacksquare$

\bigskip

The graphs with $\gamma=n/2$ are known; they are coronas or unions of
4-cycles, see \cite{HHS}. If the corona of $H$ is triangle-free, then so is
$H$. Thus a connected triangle-free graph $G$ such that $\Delta(G)\leq3$
satisfies $\gamma(G)=\gamma^{\infty}(G)$ if and only if $G=C_{4}$, or $G$ is
the corona of $P_{n}$, $n\geq1$, or of $C_{n}$, $n\geq4$. We improve this
result for triangle-free graphs. Again we need a lemma about the existence of
an EDS in which every vertex has an external private neighbor.

\begin{lemma}
\label{LemT-free}If $G$ is a triangle-free graph without isolated vertices
such that $\gamma^{\infty}(G)=\gamma(G)$, then $G$ has a minimum EDS $D$ such
that $\operatorname{epn}(v,D)\neq\varnothing$ for each $v\in D$.
\end{lemma}

\noindent\textbf{Proof}.\hspace{0.1in} Let $D$ be a minimum EDS of $G$ that
maximizes the number of edges in $G[D]$. Suppose $\operatorname{epn}%
(u,D)=\varnothing$ for some $u\in D$. Since $D$ is a minimum dominating set,
$u$ is isolated in $G[D]$. Since $\deg u\geq1$, $u$ is adjacent to a shared
vertex $w$. If $u$ defends $w$, then $D^{\prime}=(D-\{u\})\cup\{w\}$ is an EDS
such that $G[D^{\prime}]$ has more edges than $G[D]$, a contradiction.
Therefore $w$ is defended by $v\in D$ such that $\operatorname{epn}%
(v,D)\neq\varnothing$. By Fact \ref{FactEDS}, $G[\{v,w\}\cup\operatorname{epn}%
(v,D)]\cong K_{n}$ for some $n\geq3$, which is impossible in a triangle-free
graph. \hfill$\blacksquare$

\begin{theorem}
Let $G$ be a triangle-free graph with $\gamma^{\infty}(G)=\gamma(G)$. Then
$\gamma^{\infty}(G)=\theta(G)$.
\end{theorem}

\noindent\textbf{Proof}.\hspace{0.1in} Assume without loss of generality that
$G$ has no isolated vertices and let $D$ be a minimum EDS such that
$\operatorname{epn}(v,D)\neq\varnothing$ for each $v\in D$; such a set $D$
exists by Lemma \ref{LemT-free}. By Fact \ref{FactEDS}, $\{v\}\cup
\operatorname{epn}(v,D)$ forms a clique. Since $G$ is triangle-free,
$|\operatorname{epn}(v,D)|=1$ for each $v\in D$; say $\operatorname{epn}%
(v,D)=\{v^{\prime}\}$. Let $C$ be the set of all shared vertices. If
$C=\varnothing$, then we are done, so assume $C\neq\varnothing$; say $w\in C$.
Since $D$ is an EDS, $w$ is defended by some vertex $v\in D$. But then Fact
\ref{FactEDS} implies that $w$ is adjacent to $v^{\prime}$, that is,
$\{v,v^{\prime},w\}$ forms a triangle, a contradiction. \hfill$\blacksquare$

\section{Open Problems}

\label{Sec_Open}We consider Questions \ref{main1} and \ref{main2} to be
fundamental questions in the study of eternal domination.

\begin{question}
\label{main1} Does there exist a graph $G$ such that $\gamma(G)=\gamma
^{\infty}(G)$ and $\gamma(G)<\theta(G)$?
\end{question}

\begin{question}
\label{main2} Does there exist a triangle-free graph $G$ such that
$\gamma^{\infty}(G)=\alpha(G)<\theta(G)$?
\end{question}

We do not know of similar questions to Questions \ref{main1} and \ref{main2}
in the $\mathrm{m}$-eternal domination problem. For example, $\gamma
(C_{n})=\gamma_{\mathrm{m}}^{\infty}(C_{n})=\alpha(C_{n})<\theta(C_{n})$ when
$n\in\{5,7\}$ (and, of course, $C_{n}$ is triangle-free for $n>3$).

There exist triangle-free graphs $G$ with $\theta(G)=\gamma^{\infty}(G)$ and
$\alpha(G)<\theta(G)$; $C_{5}$ is one example. Infinitely many graphs that are
not triangle-free with the property that $\alpha(G)=\gamma^{\infty}%
(G)<\theta(G)$ are described in \cite{KM2}, as well as graphs with
$\alpha(G)<\gamma^{\infty}(G)<\theta(G)$. It remains open to characterize all
graphs having $\gamma(G)=\gamma^{\infty}(G)$.

\begin{question}
Is it true for all planar graphs $G$ that $\gamma(G)=\gamma^{\infty}(G)$
implies $\gamma(G)=\theta(G)$?
\end{question}

Determining additional classes of graphs for which $\gamma(G)=\gamma
_{\mathrm{m}}^{\infty}(G)$, $\gamma_{\mathrm{m}}^{\infty}(G)=\alpha(G)$, or
$\gamma_{\mathrm{m}}^{\infty}(G)=\theta(G)$ is also an interesting direction
for future work. As mentioned in Section \ref{Sec_Small_alpha}, if
$\alpha(G)=2$, then $\gamma_{\mathrm{m}}^{\infty}(G)=2$, and if $\alpha(G)=3$,
then $\gamma_{\mathrm{m}}^{\infty}(G)\in\{2,3\}$. Proposition
\ref{m-two-converse} gives a sufficient condition for $\gamma_{\mathrm{m}%
}^{\infty}(G)$ to equal $2$ while $\alpha(G)=3$. The following problem could
be a starting point for an investigation into graphs that satisfy
$\gamma_{\mathrm{m}}^{\infty}(G)=\alpha(G)$.

\begin{problem}
Characterize the class of graphs $G$ such that $2=\gamma_{\mathrm{m}}^{\infty
}(G)<\alpha(G)=3$ (equivalently $\gamma_{\mathrm{m}}^{\infty}(G)=\alpha(G)=3$).
\end{problem}

Sixty one Cayley graphs of nonabelian groups for which $\gamma_{\mathrm{m}%
}^{\infty}(G)=\gamma(G)+1$ were discovered by Braga et al.~in \cite{BdSL}.
Disjoint unions of these graphs give examples of Cayley graphs for which the
difference $\gamma_{\mathrm{m}}^{\infty}(G)-\gamma(G)$ can be an arbitrary
positive integer, but at present there is no similar result for connected
Cayley graphs.

\begin{question}
Does there exist a connected Cayley graph $G$ such that $\gamma_{\mathrm{m}%
}^{\infty}(G)>\gamma(G)+1$? Can the difference $\gamma_{\mathrm{m}}^{\infty
}(G)-\gamma(G)$ be arbitrary for connected Cayley graphs?
\end{question}

\begin{problem}
Find an infinite class of connected Cayley graphs such that $\gamma
_{\mathrm{m}}^{\infty}(G)>\gamma(G)$.
\end{problem}

The next question relates to Fact \ref{Fact_Necessary} and Proposition
\ref{Prop_gamma2}.

\begin{question}
For $k\geq3$, which graphs $G$ satisfy $\gamma^{\infty}(G)=k$, or
$\gamma_{\mathrm{m}}^{\infty}(G)=k$, if and only if every vertex of $G$ is in
a dominating set of size $k$?
\end{question}

Let $G\boksie H$ denote the Cartesian product of $G$ and $H$. An interesting
conjecture is that of Finbow and Klostermeyer \cite{KMSurvey}, who conjectured
there exists a constant $c$ such that $\gamma_{\mathrm{m}}^{\infty}%
(P_{n}\boksie P_{n})\leq\gamma(P_{n}\boksie
P_{n})+c$, for all $n$. We state another conjecture.

\begin{conjecture}
\label{c1} Let $G$ be a graph such that $\theta(G)=\gamma^{\infty}(G)$. Then
$\theta(G \boksie K_{2})=\gamma^{\infty}(G \boksie K_{2})$.
\end{conjecture}

Perhaps Conjecture \ref{c1} is also true if $K_{2}$ is replaced with any tree.
Similar statements for $\gamma_{\mathrm{m}}^{\infty}(G)$ do not seem to be
true. For example, let $G$ be a graph such that $\gamma(G)=\gamma_{\mathrm{m}%
}^{\infty}(G)$. In many cases, $\gamma(G\boksie K_{2})=\gamma_{\mathrm{m}%
}^{\infty}(G\boksie K_{2})$. But $\gamma(K_{2,3}-e\boksie K_{2})=3<\gamma
_{\mathrm{m}}^{\infty}(K_{2,3}-e\boksie K_{2})=4.$ Likewise, if we replace
$\gamma$ with $\theta$ in this, we find the following example: $\theta
(C_{4}\boksie K_{2})=4>\gamma_{\mathrm{m}}^{\infty}(C_{4}\boksie K_{2})=3$.

One might consider Vizing-like conjectures by asking whether $\gamma
_{\mathrm{m}}^{\infty}(G\boksie H)\geq\gamma_{\mathrm{m}}^{\infty}%
(G)\ast\gamma_{\mathrm{m}}^{\infty}(H)$, for all $G,H$. But this is not true
in general, as $\gamma_{\mathrm{m}}^{\infty}(P_{3}\boksie P_{3})=3<\gamma
_{\mathrm{m}}^{\infty}(P_{3})\ast\gamma_{\mathrm{m}}^{\infty}(P_{3})=4$. A
proof that $\gamma_{\mathrm{m}}^{\infty}(P_{3}\boksie P_{3})=3$ can be found
in \cite{GK2}. Perhaps $\gamma_{\mathrm{m}}^{\infty}(G\boksie H)\geq
\max\{\gamma_{\mathrm{m}}^{\infty}(G)\ast\gamma(H),\gamma(G)\ast
\gamma_{\mathrm{m}}^{\infty}(H)\}$, for all $G,H$?


However, the Vizing-like problem for eternal domination seems challenging.

\begin{question}
Is it true for all graphs $G,H$ that $\gamma^{\infty}(G\boksie H)\geq
\gamma^{\infty}(G)\ast\gamma^{\infty}(H)$?
\end{question}


\noindent{\large \textbf{{Acknowledgements}}}

The authors wish to thank Gary MacGillivray for some helpful discussions, and
the referees for valuable comments and suggestions.

\end{document}